\documentclass[10pt]{article}
\usepackage{epsfig,amsmath,amsthm,amssymb,verbatim,fullpage}
\theoremstyle{plain}
\newtheorem{theorem}{Theorem}[section]  
\newtheorem*{corollary}{Corollary}
\newtheorem*{oldtheorem}{Theorem}       
\newtheorem{lemma}{Lemma}[section]
\newtheorem{definition}{Definition}[section]
\newtheorem{remark}{Remark}[section]
\newtheorem{example}{Example}[section]
\def\vect#1{\overrightarrow{#1}}
\let\a\alpha  \let\b\beta  \let\c\gamma

\begin{document}

\title{How to sum up triangles}
\author{Bakharev F. \and Kokhas K.\and Petrov F.}
\date{\small June 2001}
\maketitle

\begin{abstract}
We prove configuration theorems that generalize
the Desargues, Pascal, and Pappus theorems. 
Our generalization of the Desargues theorem allows 
us to introduce the structure of an Abelian group on 
the (properly extended) set of triangles which are 
perspective from a point. In barycentric coordinates,
the corresponding group operation becomes 
the addition in~$\mathbb R^3$.\\ [3mm]


\end{abstract}

\section{Configuration theorems}

\subsection{The Desargues theorem and its generalization}

In this paper we deal with points and (straight) lines on 
the (real projective) plane. By a \emph{configuration} 
we understand a set of points and lines on this plane.

\begin{definition}
If there is a correspondence between two configurations
such that all the lines passing through the corresponding 
points meet in a point $O$, then we say that these 
configurations are \emph{perspective from the point~$O$} 
and we call $O$ \emph{the perspective center}.

If there is a correspondence between two configurations
such that all the intersection points of the corresponding 
lines lie on a line~$\ell$, then we say that these 
configurations are \emph{perspective from the line~$\ell$} 
and we call $\ell$ \emph{the perspective axis}.
\end{definition}

\begin{oldtheorem}[Desargues]
If two triangles are perspective from a point, 
they are also perspective from a line.
\end{oldtheorem}

The converse (also referred to as the ``dual'') 
Desargues theorem holds as well: if two triangles are 
perspective from a line, then they are perspective from a point.

An interested reader can find a proof of the 
Desargues theorem in the books \cite{Coxeter}, 
\cite{Kon-Fossen}, \cite{Efimov}.

Thus, the Desargues theorem states that the intersection 
points of the corresponding sides of two perspective 
triangles lie on a line. Can we say anything about the 
intersection points of \emph{non-corresponding} sides?
It turns out that we can!

\begin{figure}
\begin{center}
\epsfig{file=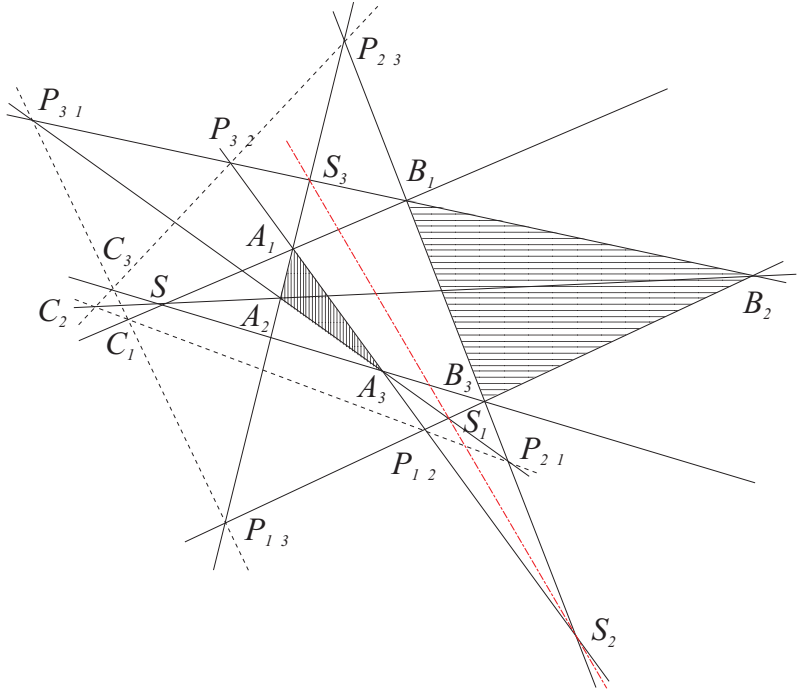,width=.7\hsize}
\caption{Desargues theorem and its generalization}
\label{fig:desarg_g}
\end{center}
\end{figure}

\begin{definition}[The main construction]
Let $A_1A_2A_3$ and $B_1B_2B_3$ be triangles 
perspective from a point~$S$. Let $(i,j,k)$ stand for a 
permutation of the numbers 1, 2,~3. Let $P_{ij}$ be the 
intersection point of the lines $A_iA_k$ and $B_jB_k$.
We denote $S_k=A_iA_j\cap B_iB_j$ and
$C_k=P_{ik}P_{ki}\cap P_{jk}P_{kj}$ 
(see~fig.~\ref{fig:desarg_g}).
\end{definition}

\begin{theorem}[Generalized Desargues theorem]
\label{thm:gen_Desarg}
The triangle $C_1C_2C_3$ is perspective to the triangles 
$A_1A_2A_3$ and $B_1B_2B_3$ from the point~$S$.
\end{theorem}

We will give two proofs of this theorem.

\begin{proof}[Proof 1]
Let us check that the triangles $\Delta_1=S_1P_{21}P_{12}$ 
and $\Delta_2=S_3P_{23}P_{32}$ are perspective from
the point~$S_2$. For this purpose we will show 
that the lines $P_{21}P_{23}$, $P_{12}P_{32}$, and $S_1S_3$
pass through the point~$S_2$. Observe that $P_{21}P_{23}$ 
is actually the line $B_3B_1$ and that $P_{12}P_{32}$ is 
actually the line $A_1A_3$. Therefore these two lines 
intersect in the point~$S_2$. It remains to notice that 
$S_1$, $S_2$, and $S_3$ are collinear by an application 
of the Desargues theorem to the 
triangles $A_1A_2A_3$ and $B_1B_2B_3$.

Thus, the triangles $\Delta_1$ and $\Delta_2$ are 
perspective from a point. Hence, by the Desargues theorem,
the intersection points of the corresponding sides,
\begin{align*}
A_2 &= S_1P_{21}\cap S_3P_{23},\\
B_2 &= S_1P_{12}\cap S_3P_{32},\\
C_2 &= P_{21}P_{12}\cap P_{23}P_{32}, 
\end{align*}
are collinear. Analogously, $C_1\in A_1B_1$, $C_3\in A_3B_3$, 
which completes the proof. 
\end{proof}

\begin{proof}[Proof 2]
We will need the following fact about algebraic curves
(see~\cite{Solovjev}, chapter~1, \S~1, Theorem~4 for $n=4$):
\begin{quote}
Let $A_{ij}$ be (pairwise distinct) intersection points 
of lines $p_i$, $q_j$, $1\leqslant i, j\leqslant 4$.
If all the following points 
$$
\begin{array}{cccc}
A_{11},& A_{12},& A_{13},& A_{14}, \\
A_{21},& A_{22},& A_{23},& A_{24}, \\
A_{31},& A_{32},& A_{33}, \\
A_{41},& A_{42} \\
\end{array}
$$
belong to an algebraic curve of degree four, then the 
remaining points $A_{34}$, $A_{43}$, $A_{44}$ belong to 
this curve as well. 
\end{quote}
Consider the triangles $A_1A_2A_3$ and $B_1B_2B_3$ that
are perspective from the point~$S$. Construct the points
$P_{ij}$. Define the points $C_i$ not as in the main
configuration but as follows:
$C_1'= P_{12}P_{21}\cap A_1B_1$,
$C_2'= P_{12}P_{21}\cap P_{23}P_{32}$,
$C_3'= P_{23}P_{32}\cap A_3B_3$.
Let us denote straight lines in an over-determined way
by listing those points of the obtained configuration
which belong to these lines. Consider the lines 
\begin{align*}
p_1&=P_{23}A_1A_2P_{13},  &q_1&=P_{31}A_2A_3P_{21},  
&r_1&=P_{31}P_{13},      \\
p_2&=P_{31}P_{32}B_1B_2,  &q_2&=P_{13}P_{12}B_3B_2,  
&r_2&=P_{32}A_1A_3P_{12},\\
p_3&=C_3' S A_3B_3,       &q_3&=C_1'SA_1B_1,         
&r_3&=P_{23}B_1B_3P_{21},\\
p_4&=C_2'C_1'P_{12}P_{21},&q_4&=C_2'C_3'P_{32}P_{23},
&r_4&=SA_2B_2.
\end{align*}
We see that the intersection points of $p_i$ and $q_i$
give 16 points of our configuration. All of them,
except maybe $C_1'$, $C_2'$, and $C_3'$, belong to
a curve of degree four that is the union 
of the four lines
$r_1$, $r_2$, $r_3$, and $r_4$. In the notations of the
above mentioned theorem from \cite{Solovjev} we have  
$C_1'=A_{43}$, $C_2'=A_{44}$, $C_3'=A_{34}$. Therefore,
these points must also belong to the union 
of the lines~$r_i$. 
It is easy to check that none of the lines~$r_i$ contains
five points of our configuration. Hence the points $C_3'$
and $C_1'$ lie on $r_1$; and the point $C_2'$ lies
on~$r_4$. This completes the proof.
\end{proof}

In the book \cite{Kon-Fossen} (\S~22), the Reye 
configuration with parameters (12, 4, 16, 3) is described. 
The vertices of the configuration are the vertices of 
a cube, its center, and three points at infinity 
which correspond to the directions of the edges of the 
cube. The lines of the configuration are the edges and 
the four main diagonals of the cube.
It is not difficult to verify that the configuration
described in Theorem~\ref{thm:gen_Desarg} is
dual to the Reye configuration.

\begin{quote}
\small
Remark. The Pascal theorem can be proven by
similar considerations. In fact, it is a theorem on
associativity of the addition operation for points
on an elliptic curve. Unfortunately our theorem~%
\ref{thm:gen_Desarg} does not seem to be associated 
with any operation for objects related to a curve of 
degree four. We do not have even a configuration 
theorem for a generic curve of degree four.
\end{quote}

\subsection{The Pappus and Pascal theorems} 

In this subsection we will prove generalizations
of the Pappus and Pascal theorems. Readers
familiar with the polar transform and dual theorems 
can easily generalize the Brianchon theorem along 
the same lines.

\begin{oldtheorem}[Pascal]
If a simple hexagon is inscribed in a conic section,
then the three intersection points of opposite sides 
of the hexagon are collinear. 
\end{oldtheorem}

\begin{oldtheorem}[Pappus]
Let the points $A_1$, $A_2$, $A_3$ lie on one
line and the points $B_1$, $B_2$, $B_3$ lie on 
another  line. Then the intersection points of the 
pairs of lines $A_1B_2$ and $A_2B_1$, $A_1B_3$ and 
$A_3B_1$, $A_2B_3$ and $A_3B_2$ are collinear.
\end{oldtheorem}

Proofs of these theorems are given, e.g., in
\cite{Coxeter}, \cite{Solovjev}.

\begin{figure}
\begin{center}
\epsfig{file=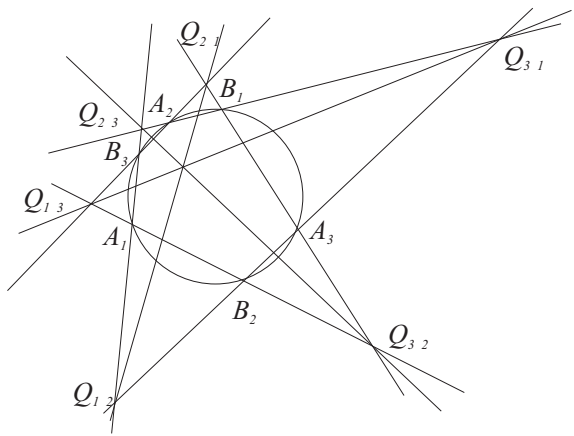,width=.5\hsize}
\caption{A generalization of the Pascal theorem}\label{pascal_g}
\end{center}
\end{figure}

We will obtain a generalization of the Pascal theorem if
we consider the intersection points of non-opposite sides.
Analogously, in order to generalize the Pappus theorem,
we consider the intersection points of ``non-symmetric''
 lines. We will use the following notations.
Let $A_1$, $A_2$, $A_3$, $B_1$, $B_2$, $B_3$ denote
the initial points (in the Pappus theorem) or the
hexagon's vertices (in the Pascal theorem; one should
pay attention to the order of the vertices, see 
fig.~\ref{pascal_g}). For every permutation $(i,j,k)$ 
of the numbers 1, 2,~3, let $Q_{ij}$ denote
the intersection point of the lines $A_iB_k$ and~$B_jA_k$.

\begin{theorem}[Generalized Pascal theorem]
Let $A_1B_3A_2B_1A_3B_2$~be an inscribed hexagon.
Then the  lines $Q_{12}Q_{21}$, $Q_{13}Q_{31}$, 
and $Q_{23}Q_{32}$ meet in a point.
\end{theorem}

\begin{theorem}[Generalized Pappus theorem]
Let the conditions of the Pappus theorem be fulfilled.
Then the  lines $Q_{12}Q_{21}$, $Q_{13}Q_{31}$, 
and $Q_{23}Q_{32}$ meet in a point.
\end{theorem}

\begin{proof}
Consider the triangles $Q_{12}Q_{23}Q_{31}$ and
$Q_{21}Q_{32}Q_{13}$. The intersection points of
the opposite sides of these triangles are exactly
the points dealt with in the Pascal theorem
(for the hexagon $A_1B_3A_2B_1A_3B_2$). Therefore 
these triangles are perspective from a line and 
hence from a point. This completes the proof.
\end{proof}

\begin{quote}\small
Let us give another, well-known  way to generalize 
the Pascal theorem. 

\begin{oldtheorem}[Another generalization of the
 Pascal theorem]
Denote $S_{ik}=A_iA_j\cap B_jB_k$. Then the
three  lines $l_i=S_{jk}S_{kj}$ meet in a point.
\end{oldtheorem}

\begin{proof}
Denote $S=A_3B_3\cap A_2B_2$. By the Pascal theorem
(using a slightly different enumeration of points), 
the points $S$, $S_{32}$, and $S_{23}$ are collinear.
Applying the Pappus theorem to 
the triples $B_3$, $S_{13}$, $B_2$ and $A_2$, 
$S_{21}$, $B_3$, we infer that the point 
$U=B_2S_{21}\cap A_3S_{13}$ belongs to~$l$. Applying 
the Pappus theorem to the triples $B_2$, $S_{12}$, 
$S_{13}$ and $A_3$, $S_{31}$, $S_{21}$, we infer that
the point $Z=S_{12}S_{21}\cap S_{13}S_{31}$ belongs
to~$l$. This completes the proof.
\end{proof}

\end{quote}

\subsection{Projective duality}

Let us consider the main construction on the 
projective plane. Notice that projective 
transformations preserve all the elements of the 
construction. Moreover, from the projective viewpoint, 
the choice of the point~$S$ and of the  lines 
$\ell_1$, $\ell_2$, and $\ell_3$ does not lead to 
loss of generality. Indeed, using a suitable
projective transformation, we can map this triple 
of concurrent  lines into any other triple of 
concurrent  lines. In particular, the structure of 
the ``additive group of triangles''
described in Section ~\ref{sec:AddGroup} 
does not depend on this choice.

\begin{figure}
\begin{center}
\epsfig{file=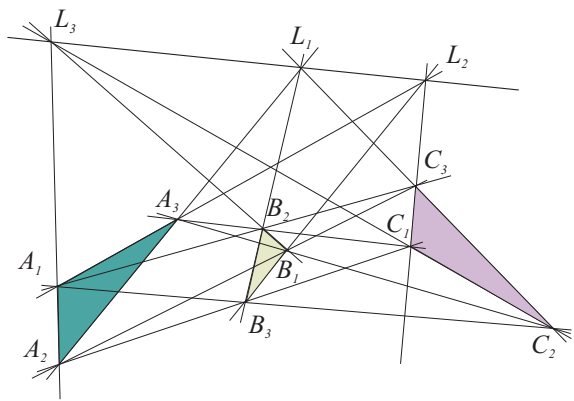,width=.5\hsize}
\caption{The dual to the main construction 
  }\label{fig:dual_desarg_g}
\end{center}
\end{figure}

Let us choose a polar transform of the projective plane.
Applying it to all the elements of the main configuration,
we obtain a theorem dual to Theorem~\ref{thm:gen_Desarg}:

\begin{theorem}\label{thm:dual_gen_Desarg}
Let the triangles $A=A_1A_2A_3$ and $B=B_1B_2B_3$ be
perspective from a line~$s$. 
For $i=1$, $2$, $3$, let $L_i$ denote the intersection 
points of the corresponding sides: $L_i=A_kA_j\cap B_kB_j$.
Let $C_k=A_iB_j\cap A_jB_i$. Then the triangle 
$C=C_1C_2C_3$ is perspective to the triangles
$A_1A_2A_3$ and $B_1B_2B_3$ from the line~$s$
(see fig.~\ref{fig:dual_desarg_g}).
\end{theorem}

Since a polar transform is involutory, 
Theorem~\ref{thm:dual_gen_Desarg} follows from
Theorem~\ref{thm:gen_Desarg}. 

\begin{definition}[Central and axis models]
The construction used in Theorem~\ref{thm:gen_Desarg}
(fig.~\ref{fig:desarg_g}) will be referred to
as the \emph{central model}. The construction used 
in Theorem~\ref{thm:dual_gen_Desarg}
(fig.~\ref{fig:dual_desarg_g}) will be referred to
as the \emph{axis model}.
\end{definition}

In the central model we fix the point~$S$ and the three 
lines $\ell_1$, $\ell_2$, $\ell_3$ that pass through
this point. In the axis model we fix the  
line~$s$ and the three points $L_1$, $L_2$, $L_3$ that
belong to this line. For brevity of notations, we 
denote by~$A$ the triangle $A_1A_2A_3$ (where 
$A_i\in\ell_i$ for the central model; $A_iA_j$ 
contains the point $L_k$ for the axis model).

\section{Additive group of triangles}
\label{sec:AddGroup}

Let us remind how the addition of points on an elliptic 
curve is defined. First, one introduces certain 
``geometric'' construction. Given two points $A$ 
and $B$ on an elliptic curve, this construction allows 
us to construct one more point of the curve, $C$. 
Namely, $C$~is the third point where
the  line $AB$ intersects the elliptic
curve. Denote this as $C=A\cdot B$. Now the addition
operation is defined as follows. We fix an arbitrary
point $E$ on the curve and put $A+B=E\cdot (A\cdot B)$. 
The point $E$ plays the role of a neutral element
(the ``zero'').

However, this method requires a discussion of certain
technical details. And it turns out that it is natural
to consider the curve not on a plane but on the 
complex projective plane. The definition of the point
$C$ requires further conventions if $A=B$ or if the
line $AB$ is tangent to the curve. Finally, one often 
chooses as $E$ a point at infinity. This
leads to a somewhat ``mysterious'' definition of the
sum: having constructed the point $A\cdot B$, one
symmetrically reflects it with respect to the axis~$x$.
The top of the theory is a construction of elliptic
functions that are homomorphisms between the additive
group of points of an elliptic curve and the complex torus.

In this section we will define an addition operation 
on the set of perspective triangles. The configuration 
theorem given in the previous subsection (Theorem~%
\ref{thm:dual_gen_Desarg}) allows us to construct a 
triangle by two given triangles. This however requires 
some technical comments, for instance, in the
case of coinciding triangles. In order to define 
the addition we need only to repeat the operation
with a fixed triangle. However, to introduce
such an operation, it requires a certain completion
of the set of triangles. It is convenient to choose
an ``infinitely remote triangle'' as the 
fixed element. Then the barycentric coordinates
will play the role of an elliptic function.

\subsection{Sum of perspective triangles}

\begin{definition}
\label{def:predsumma}
Let $A$ and $B$ be two triangles perspective from a 
point~$S$ or a  line~$s$. By the generalized Desargues 
theorem, these triangles define the third triangle, $C$, 
which we call the \emph{pre-sum} of the triangles $A$ 
and $B$. We denote this as $C=A\boxplus B$.
\end{definition}

Let us describe some obvious properties of 
this operation.

\begin{lemma}\label{lemma:quasygroup}
The operation ``pre-sum'' possesses the following 
properties:

1) $A\boxplus B=B\boxplus A$.

2) If $A\boxplus B=C$, then $A\boxplus C = B$ and
 $B\boxplus C=A$.
\end{lemma}

\begin{remark}
The set of triangles is not closed with respect to
the operation ``pre-sum''. For example, the reader 
can easily see from the fig.~\ref{fig:desarg_g}
that all the three vertices of the triangle $C$
can coincide with the point $S$ for a certain 
choice of the triangles $A$ and $B$.
\end{remark}

\begin{definition}
\label{def:summa_general}
Let $F$ be an arbitrary fixed triangle. 
The \emph{sum} of triangles $A$ and $B$ is defined 
as $A+B=F\boxplus(A\boxplus B)$ whenever 
this formula makes sense. 
\end{definition}

\begin{remark}
A set equipped with an operation that has properties
as in Lemma~\ref{lemma:quasygroup} is called a quasi-group.
Let us stress that the axioms of a quasi-group do not 
imply associativity of the operation 
$A+B=F\boxplus(A\boxplus B)$. For instance, applying
the Pappus theorem, we can construct from two triples
of collinear points $A$ and $B$ the third triple,
$C=A\boxplus B$. It is easy to verify that such an
operation satisfies Lemma~\ref{lemma:quasygroup}. 
But then a ``sum'' arising by an analogy with Definition
\ref{def:summa_general} is not associative.
\end{remark}

Further we will deal with the axis model. Let the
perspective axis, $s$, be at infinity (this
can be always achieved with the help of an appropriate 
projective transformation). For technical reasons, 
it is convenient to choose $F$ to be an infinitely 
remote triangle (this notion will be clarified later). 
Thus, in the rest of the paper we will use the 
following definition of the sum.

\begin{definition}
\label{def:summa}
Let $A$ and $B$ be triangles perspective from a line $s$, 
and let $C=A\boxplus B$. Let $C_0$ denote the mass center 
of the triangle~$C$. Let $D$ be a triangle that is centrally 
symmetric to the triangle $C$ with respect to the point~$C_0$. 
We call $D$ the \emph{sum} of the triangles $A$ and $B$ and 
write this as $D=A+B$.
\end{definition}

We will prove below that this operation is commutative
and associative (which justifies calling it a sum). 
We will show that the symmetry with respect to the 
mass center is actually a computation of the pre-sum of
a given triangle with a certain ``infinitely remote
triangle'' (see Example~\ref{ex:a_plus_000}). The latter 
will play the role of the neutral element (the ``zero'') 
for the addition operation. Let us remark that the operation 
defined in Definition~\ref{def:summa_general} is also 
commutative and associative and gives rise to the same 
Abelian group. The proof of this is completely 
analogous to that in the case of the sum defined 
in Definition~\ref{def:summa}.

\subsection{Barycentric coordinates and 
 the space of triangles}

\begin{figure}
\begin{center}
\epsfig{file=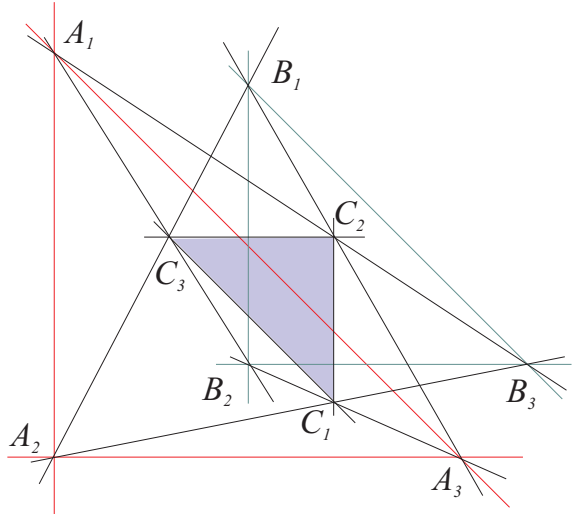,width=.5\hsize}
\caption{The dual construction with an axis at infinity}
\label{fig:dual_desarg_inf}
\end{center}
\end{figure}

Applying an appropriate projective transformation,
we can always achieve that the  line $s$ in 
the axis model is at infinity. In this case
the corresponding sides of the triangles $A$, $B$, 
and $C=A\boxplus B$ are parallel 
(see fig.~\ref{fig:dual_desarg_inf}). Below we
understand by {\sf triangles} only triangles which
have the same directions of sides.

Choose any {\sf triangle}~$E=E_1E_2E_3$. Below we
understand by coordinates of a point on the affine
plane the barycentric coordinates of this point
with respect to the {\sf triangle}~$E$ (see~\cite{Balk}).
With the help of these coordinates we will introduce
coordinates on the set of {\sf triangles}.

\begin{definition}[Barycentric coordinates on the set 
 of {\sf triangles}]
\label{def:BarCoordinates}
Every {\sf triangle}~$D$ is uniquely determined by the
coordinates $(d_1,d_2,d_3)$ (where $d_1+d_2+d_3=1$) of 
its mass center and by the coefficient~$d$ of the
only homothecy which transforms $D$ into $E$ (the parallel
translation is by definition a homothecy with $d=1$).
Consequently, this {\sf triangle} is uniquely determined
by a triple of \emph{ parameters},
$(\delta_1,\delta_2,\delta_3)$, such that
$$
\delta_1:\delta_2:\delta_3=d_1:d_2:d_3\,, \qquad
\delta_1+\delta_2+\delta_3=d\,.
$$
The parameters $(\delta_1,\delta_2,\delta_3)$ are 
called the \emph{barycentric coordinates} of the 
{\sf triangle}~$D$.
\end{definition}

The connection between the coordinates 
($d_1$, $d_2$, $d_3$;~$d$) and the barycentric coordinates
$(\delta_1,\delta_2,\delta_3)$ is given by
$$
\delta_i=d\cdot d_i \qquad\
d=\delta_1+\delta_2+\delta_3; \quad
d_i=\frac{\delta_i}d= 
\frac{\delta_i}{\sum \delta_i}\quad (i=1,2,3).
$$

\begin{example}
\label{ex:pereschet_coordinat}
Let us find the coordinates of the vertices of a 
{\sf triangle} $D=D_1D_2D_3$ with the mass center 
$D_0$ in terms of its barycentric coordinates
$\delta_1$, $\delta_2$, $\delta_3$. For this purpose
we notice that the vector $\overrightarrow{D_0D_3}$
is homothetic to the vector $\overrightarrow{E_0E_3}$ 
with a coefficient~$\frac1d$:
$$
D_3=D_0+\frac 1d\,(E_3-E_0)=\frac1d\, 
 (\delta_1,\delta_2,\delta_3)+
\frac1d\, \Bigr(\!-\frac{1}3,-\frac{1}3,\frac23\Bigl)=
\Bigl(\frac{\delta_1-1/3}d, \frac{\delta_2-1/3}d, 
 \frac{\delta_3+2/3}d\Bigr).
$$
Coordinates of the other vertices are computed
analogously.
\end{example}

Let us remark that there exists no (geometric) 
{\sf triangle} if the sum of the coordinates
vanishes. 

\begin{definition} \label{def:TriangleSpace}
Let us extend the set of geometric {\sf triangles} by 
adding formal elements whose barycentric coordinates
are given by parameters $(\delta_1, \delta_2, \delta_3)$
with $\delta_1+\delta_2+\delta_3=0$. We call the 
resulting space the \emph{space of triangles}; it
is formally isomorphic to $\mathbb R^3$.
\end{definition}

Triangles which have small sum of their coordinates
are of a ``big size''. It is therefore natural to
expect that the formal elements admit an interpretation
as ``infinitely remote'' triangles
(see Section~\ref{sec:geom_details}).

\subsection{The main theorem}

\begin{theorem}\label{thm:associativity}
The operation ``$+$'' defined in
Definition~\ref{def:summa} coincides, in the barycentric
coordinates, with the addition operation in~$\mathbb R^3$.
\end{theorem}

\begin{proof}
Let $A$ and $B$ be two {\sf triangles} with 
coordinates $(\a_1,\a_2,\a_3)$, $\sum\a_i=a\ne0$, 
and $(\b_1,\b_2,\b_3)$, $\sum b_i=b\ne0$,
respectively. Let their pre-sum $C=A\boxplus B$
be also a {\sf triangle}, i.e., $\sum\c_i=c\ne0$ 
for coordinates $(\c_1,\c_2,\c_3)$ of $C$.
Let us find the coordinates of the {\sf triangle}~$C$.
The triangles $A_1C_3A_2$ and $B_2C_3B_1$ 
(see fig.~\ref{fig:dual_desarg_inf}) are homothetic;
the homothecy coefficient is  $A_2A_1/B_2B_1=\frac{1/a}{1/b}$.
Hence $a\cdot A_1C_3=b\cdot C_3B_2$. This
implies that the point $C_3$ is the mass center 
of the points $(a,A_1)$ and $(b,B_2)$. That is,
\begin{multline}
\label{koord_predsummy}
C_3=\frac1{a+b}\,(a\cdot A_1+b\cdot B_2)=
  \frac1{a+b}\,\bigl((\a_1+2/3, \a_2-1/3, \a_3-1/3)+
  (\b_1-1/3, \b_2+2/3, \b_3-1/3) \bigr)=\\
  =\Bigl(\frac{-(\a_1+\b_1)-1/3}{-(a+b)}, 
  \frac{-(\a_2+\b_2)-1/3}{-(a+b)},
  \frac{-(\a_3+\b_3)+2/3}{-(a+b)}\Bigr)\,.
\end{multline}
Coordinates of the points $C_1$ and $C_2$ are found
in the same way. Notice that
$$
C_2-C_1=\frac1{-(a+b)}\,(-1,1,0)=\frac1{-(a+b)}\,(E_2-E_1)\,.
$$
Therefore, the sum of coordinates of the {\sf triangle}~$C$
(which is the homothecy coefficient of $C\rightarrow E$) 
equals to $-(a+b)$. Multiplying this quantity by the 
coordinates of the mass center $C_0$ of the 
{\sf triangle}~$C$, we will obtain the coordinates 
of~$C$. The coordinates of~$C_0$ are given by
$$
\frac{C_1+C_2+C_3}{3}=\left(
\frac{-(\a_1+\b_1)}{-(a+b)},
\frac{-(\a_2+\b_2)}{-(a+b)},
\frac{-(\a_3+\b_3)}{-(a+b)}
\right)\,.
$$
Therefore, the coordinates of the triangle $C$ are
\begin{equation}\label{eqn:coord_predsummy}
(-(\a_1+\b_1), -(\a_2+\b_2), -(\a_3+\b_3))\,.
\end{equation}

It is clear that the central symmetry transformation 
with respect to the mass center reverses the coordinates 
of a {\sf triangle}. Thus, coordinates of the sum
$A+B$ equal to
\begin{equation}\label{star}
(\a_1+\b_1, \a_2+\b_2, \a_3+\b_3) \,.
\end{equation}

\end{proof}

\begin{corollary}
The operation ``$+$'' defined in 
Definition~\ref{def:summa} extends to a 
commutative and associative operation on the 
space of triangles. The resulting 
``additive group of triangles''
is isomorphic to the additive group~$\mathbb R^3$.
\end{corollary}

\section{Technical details}

\subsection{Geometric interpretation of formal triangles}
\label{sec:geom_details}

As we show below, almost all formal triangles 
admit a geometric interpretation.

\begin{definition}[Barycentric coordinates on 
 the line at infinity]
\label{def:InfBarCoordinates}
Let $X$ be a point with barycentric coordinates
$(x_1,x_2,x_3)$, $\sum x_i=1$. Then 
$\vect {OX}=\sum x_i \vect {OE_i}$ for any 
point~$O$. Similarly, for a triple of numbers 
$(y_1,y_2,y_3)$ with vanishing sum, we can 
consider vector $\sum y_i \vect {OE_i}$. 
Apparently, this vector does not depend on~$O$.
Therefore, the triple $(y_1,y_2,y_3)$ determines
a direction (in fact, a vector; the direction is
obtained after a factorization over the scalar multiplier). 
We call the homogeneous triple
$(y_1,y_2,y_3)$ the \emph{barycentric coordinates} 
of this direction (= coordinates of a point on
the line at infinity).
\end{definition}

\begin{definition}
A \emph{pseudo-triangle} with coordinates 
$(p_1,p_2,p_3)$, $\sum p_i=0$, is an ordered triple 
of directions with the following coordinates
$$
(p_1+\tfrac23, p_2-\tfrac13, p_3-\tfrac13), \quad
(p_1-\tfrac13, p_2+\tfrac23, p_3-\tfrac13), \quad
(p_1-\tfrac13, p_2-\tfrac13, p_3+\tfrac23)  \,.
$$
The points of the line at infinity
corresponding to these directions (or the directions
themselves) are called the \emph{vertices} of the 
pseudo-triangle. A formal element of the space of
triangles with coordinates $(p_1,p_2,p_3)$, 
$\sum p_i=0$, is regarded as a pseudo-triangle 
$(p_1,p_2,p_3)$.
\end{definition}

\begin{figure}
\begin{center}
\epsfig{file=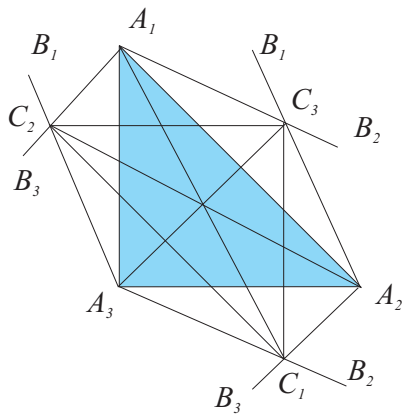,width=.35\hsize}
\caption{Addition of a {\sf triangle} with the  
 zero triangle}\label{fig:a_boxplus_000}
\end{center}
\end{figure}

\begin{example}
\label{ex:predumma_symm}
Let $A$ and $C$ be {\sf triangles} that are centrally 
symmetric to each other with respect to a point which
is not the middle of any of the sides of these 
{\sf triangles}. Let us find their pre-sum. Constructing
the main configuration, we observe that the 
lines $A_iC_j$ and $A_jC_i$ are parallel for all pairs
$i$, $j$ (see fig.~\ref{fig:a_boxplus_000}). Thus, the
``triangle'' $A\boxplus C$ is a triple of points at
infinity, i.e., it is a pseudo-triangle.
\end{example}

\begin{example}
\label{ex:000vertices}
What is a pseudo-triangle which has coordinates 
$(0,0,0)$ ? Its vertices are the directions 
$(\frac23,-\frac13,-\frac13)$, 
$(-\frac13, \frac23, -\frac13)$,
$(-\frac13, -\frac13, \frac23)$. 
Choosing the mass center of the {\sf triangle}~$E$
as the point~$O$ in 
Definition~\ref{def:InfBarCoordinates}, 
we see that these directions are the directions of
the medians of the {\sf triangle}~$E$ (or of any
other {\sf triangle}).
\end{example}

On the projective plane, the main construction 
without any additional conventions allows us to find 
the sum of a {\sf triangle} and a pseudo-triangle.

\begin{example}
\label{ex:a_plus_000}
Let us find the pre-sum of a {\sf triangle}~$A$ and 
the pseudo-triangle~$B$ with coordinates $(0,0,0)$. 
Vertices of the pseudo-triangle $B(0,0,0)$ are
the directions of the medians of the 
{\sf triangle}~$A$ (see Example~\ref{ex:000vertices}).
Constructing the main configuration, we obtain that
the pre-sum in question is simply the reflection of
the {\sf triangle}~$A$ with respect to its mass 
center (see fig.~\ref{fig:a_boxplus_000}).
\end{example}

Let us remark that, although the barycentric 
coordinates are defined up to multiplication
by a constant, pseudo-triangles with distinct
coordinates correspond to different triples of
directions.

The definition of a pseudo-triangle does not make
sense only for the three sets of parameters:
$(-\frac{2}{3},\frac{1}{3},\frac{1}{3})$,
$(\frac{1}{3},-\frac{2}{3},\frac{1}{3})$, and
$(\frac{1}{3},\frac{1}{3},-\frac{2}{3})$
since, by virtue of 
Definition~\ref{def:InfBarCoordinates},
points of the line at infinity
can not have three vanishing coordinates.

\begin{definition}
The three formal elements of the space of triangles,
which have the coordinates
$(-\frac{2}{3},\frac{1}{3},\frac{1}{3})$,
$(\frac{1}{3},-\frac{2}{3},\frac{1}{3})$, or
$(\frac{1}{3},\frac{1}{3},-\frac{2}{3})$,
are called the {\sf completely-pseudo-triangles}.
\end{definition}

Thus, the formal elements in 
Definition~\ref{def:TriangleSpace} are either 
pseudo-triangles or completely-pseudo-triangles.
This classification is exhausting.

\subsection{Corrections to the main construction} 

As before, we deal with the axis model. It is 
easy to see that certain ambiguities in the main 
construction for {\sf triangles} $A$ and $B$ 
appear in the three following cases:

1) if two sides of the {\sf triangles} $A$ and $B$ 
belong to the same line (e.g., if the points 
$A_1$, $A_2$, $B_1$, $B_2$ are collinear, then the 
lines $A_1B_2$ and  $A_2B_1$ coincide and the vertex 
$C_3$ is not defined);

2) if $A=B$ (which is a stronger degeneration
of the previous case; here none of the vertices
of the {\sf triangle}~$C$ is defined);

3) if the {\sf triangles} $A$ and $B$ are 
centrally symmetric to each other with respect
to the middle of their (common) side (say, if
$A_1A_2=B_2B_1$, then the  lines
$A_1B_2$ and $A_2B_1$ are not defined).

In the first two cases the pre-sum is a 
{\sf triangle}; in the third case it is
a completely-pseudo-triangle. Below we explain
how to modify the main construction so that
the pre-sum for these degenerate cases can
be found geometrically.

\begin{figure}
\begin{center}
\epsfig{file=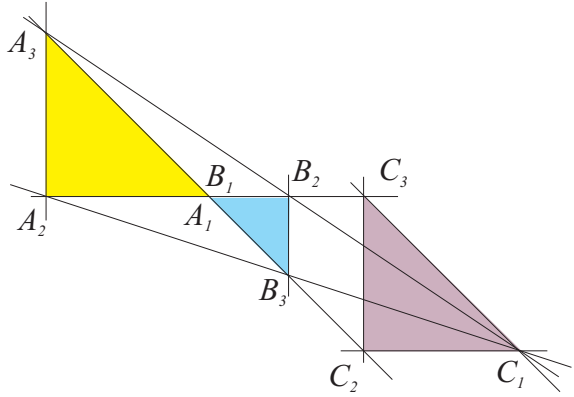,width=.5\hsize}
\caption{The simplest degenerate case of 
 pre-summation}
\label{fig:degenerate_summa}
\end{center}
\end{figure}

\subsubsection{The case of two coinciding sides}

Assume that $A_1=B_1$ but $A_2\ne B_2$ and that 
the {\sf triangles} $A$ and $B$ are not centrally
symmetric to each other 
(see fig.~\ref{fig:degenerate_summa}). Then
$C_1$ is determined by the main construction.
About $C_2$ and $C_3$ we know only that
$C_2\in A_1B_3=B_1A_3$ and $C_3\in A_1B_2=A_2B_1$. 
This information allows us to construct $C_2$ and
$C_3$ geometrically (say, 
$C_2=A_1B_3\cap C_1S_3$, where $S_3$ is the 
infinitely remote point of the  line 
$A_1A_2$). The case of $A_2=B_1$ and the case
when the points $A_1$, $B_1$, $A_2$, $B_2$ lie
on a  line can be treated analogously.
In all these cases the coordinates of the pre-sum 
are given by the formula~(\ref{eqn:coord_predsummy}).
This is clear from taking a limit.

\subsubsection{Geometric definition of  
 $A\boxplus A$}
\label{sec:a_plus_a}

\begin{figure}[cb]
\begin{center}
\epsfig{file=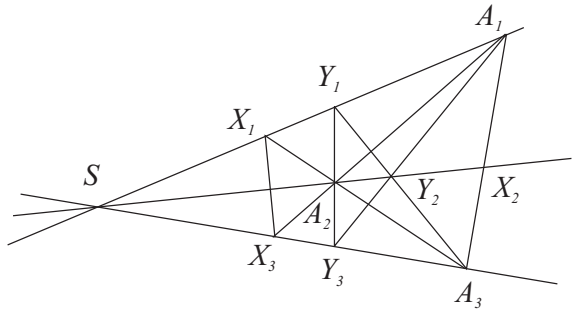,width=.5\hsize}
\caption{Adding a triangle to itself}\label{fig:a_boxplus_a}
\end{center}
\end{figure}

Let us require that the property
$A\boxplus (A\boxplus A)=A$ hold. In the central
model, this implies that the vertices of the
triangle~$A$ belong to the sides of the triangle
$A\boxplus A$. Let us construct $A\boxplus A$.
Denote $X_i=A_jA_k\cap SA_i$ (see
fig.~\ref{fig:a_boxplus_a}). Let $Y_j$ be the 
harmonic complement of the point~$S$ to the pair
$\{A_j,X_j\}$. Then three quadruples of
the form $\Gamma_j=\{S,Y_j\}\{A_j,X_j\}$ are
harmonic. Projecting the quadruple $\Gamma_j$
{}from the point $A_i$ onto the line $SA_k$, we 
obtain a harmonic quadruple $\{S,Y_k'\}\{X_k,A_k\}$,
where $Y_k'=A_iY_j\cap SA_k$. This implies that
the points $Y_k$ and~$Y_k'$ coincide since they
both are harmonic complements of the point~$S$
to the pair $\{A_k,X_k\}$. Consequently, the points 
$A_i$, $Y_j$, $Y_k$ are collinear. Therefore,
$Y_1Y_2Y_3$ is the desired triangle $A\boxplus A$. 

Let us remark that $X_1X_2X_3=A/2$ in the
sense that $X_1X_2X_3\boxplus X_1X_2X_3=A$.

For the axis model, the condition ``the vertices 
of~$A$ lie on the sides of $A\boxplus A$'' 
is replaced with the condition ``the vertices 
of $A\boxplus A$ lie on the sides of~$A$''.
This means that $A\boxplus A$ is the 
{\sf triangle} whose vertices are the middles of the sides
of~$A$. It is easy to see that this
synthetic  definition is consistent with
Definition~\ref{def:BarCoordinates} and
the formula~(\ref{eqn:coord_predsummy}).

\subsubsection{Sum of triangles that are 
symmetric with respect to the middle of
a side}
\label{sec:predumma_symm}

In this case we can suggest only a formal rule.

Assume that $A_1=B_2$ and $A_2=B_1$. Then the
main construction allows us to construct two
(infinitely remote) vertices. The third vertex, 
$A_1B_2\cap A_2B_1$, is not defined. In the
case under consideration we have
$\a_1+\a_2+\a_3=-(\b_1+\b_2+\b_3)$. Analogously
to Example~\ref{ex:pereschet_coordinat}, we 
verify that coordinates of the triangles
satisfy the following relations
$$
\a_1=-\b_1 -\frac13, \quad
\a_2=-\b_2 -\frac13, \quad
\a_3=-\b_3 +\frac23\,.
$$
In this case we can define the resulting
pre-sum as the completely-pseudo-triangle
with coordinates $(-\frac13, -\frac13, \frac23)$.

Conversely, adding the completely-pseudo-triangle
$(-\frac13, -\frac13, \frac23)$ with a 
{\sf triangle}~$A$, we obtain a {\sf triangle}
that is symmetric to $A$ with respect to the
side~$A_1A_2$.

These definitions are consistent with the
properties in Lemma~\ref{lemma:quasygroup}
and perfectly agree with the 
formula~(\ref{eqn:coord_predsummy}).

\subsection{Pre-sum of two pseudo-triangles} 

It is not difficult to verify that the pre-sum
on the projective plane of two centrally symmetric
{\sf triangles} with coordinates $(\a_1,\a_2,\a_3)$ 
and $(\b_1,\b_2,\b_3)$ (see 
Example~\ref{ex:predumma_symm}) is a 
{\sf pseudo-triangle} with coordinates
$p_i=-(\a_i+\b_i)$. That is the 
formula~(\ref{eqn:coord_predsummy}) is satisfied.
One can similarly compute coordinates of the
pre-sum of a {\sf triangle} and a pseudo-triangle.

In order to give a geometric description of 
addition of pseudo-triangles, we consider the 
following parameterization of the space of pseudo- 
(including completely-pseudo-) triangles. 
We regard pseudo- and completely-pseudo-triangles
as pre-sums of a fixed {\sf triangle} (e.g., the
{\sf triangle}~$E$) with those that are symmetric
to it. We regard these symmetric {\sf triangles} 
as parameters of pseudo-triangles.

It can be shown (for instance, employing the
barycentric coordinates) that for every 
pseudo-triangle there exists a unique 
parameterizing {\sf triangle}. Pre-summation of
pseudo-triangles is described by the following
lemma:

\begin{lemma}
If {\sf triangles} $B$ and $C$ are symmetric 
to the triangle~$E$ and pseudo-triangles 
$X$ and $Y$ are defined as the pre-sums
$E\boxplus B$ and $E\boxplus C$, then 
$X\boxplus Y=E\boxplus D$, where $D_i$ are 
the middle points of the intervals $B_iC_i$.
\end{lemma}

The proof is left for an interested reader.

\vspace*{2mm}
{\bf Acknowledgment:} The authors are
grateful to A.~Bytsko for help with translation
of the manuscript.

\vspace*{3mm} \noindent
{\sc Department of mathematics and mechanics \\
     St.Petersburg State University \\
     Bibliotechnaya pl. 2 \\
     198904 St.Petersburg, Russia } \\ [2mm]

\noindent
e-mail: kostik@kk1437.spb.edu, fedor@fp5607.spb.edu

\end{document}